\documentclass{amsart}
\usepackage{amssymb,latexsym,amsmath,amsfonts,hhline}
\usepackage{pdfsync}
\input xy
\xyoption{all}

\newcommand{\cal}{\mathcal}

\makeatletter
\renewcommand{\subsection}{\@startsection{subsection}{2}{0mm}{-2mm}{-2mm}{\bf\normalsize}}

\def\sbsnt#1{\subsection{#1}}
\makeatother

\newtheorem{formula}{}[section]
\newtheorem{definition}[formula]{Definition}
\newtheorem{corollary}[formula]{Corollary}
\newtheorem{remark}[formula]{Remark}
\newtheorem{lemma}[formula]{Lemma}
\newtheorem{proposition}[formula]{Proposition}
\newtheorem{theorem}[formula]{Theorem}

\def\thrm{\begin{theorem}}
\def\thrml#1{\begin{theorem}\label{#1}}
\def\ethrm{\end{theorem}}
\def\rmrk{\begin{remark}}
\def\rmrkl#1{\begin{remark}\label{#1}}
\def\ermrk{\end{remark}}
\def\dfntn{\begin{definition}}
\def\dfntnl#1{\begin{definition}\label{#1}}
\def\edfntn{\end{definition}}
\def\nmrt{\begin{enumerate}}
\def\enmrt{\end{enumerate}}
\def\tm#1{\item[{\rm (#1)}]}
\def\qtn{\begin{equation}}
\def\qtnl#1{\begin{equation}\label{#1}}
\def\eqtn{\end{equation}}
\def\lmm{\begin{lemma}}
\def\lmml#1{\begin{lemma}\label{#1}}
\def\elmm{\end{lemma}}
\def\prpstn{\begin{proposition}}
\def\prpstnl#1{\begin{proposition}\label{#1}}
\def\eprpstn{\end{proposition}}

\def\crllr{\begin{corollary}}
\def\crllrl#1{\begin{corollary}\label{#1}}
\def\ecrllr{\end{corollary}}
\def\css{\begin{cases}}
\def\ecss{\end{cases}}
					
\def\proof{\noindent{\bf Proof}.\ }
					
\def\cA{{\cal A}}

\def\cH{{\cal H}}

\def\cS{{\cal S}}
\def\cT{{\cal T}}

\def\cX{{\cal X}}
\def\cY{{\cal Y}}

\def\fE{{\mathfrak E}}
\def\fH{{\mathfrak H}}
\def\fU{{\mathfrak U}}
\def\fL{{\mathfrak L}}
\def\fP{{\mathfrak P}}
\def\fS{{\mathfrak S}}

\DeclareMathOperator{\aut}{Aut}
\DeclareMathOperator{\alt}{Alt}

\DeclareMathOperator{\cay}{Cay}

\DeclareMathOperator{\hol}{Hol}
\DeclareMathOperator{\id}{id}
\DeclareMathOperator{\im}{im}
\DeclareMathOperator{\Inn}{Inn}
\DeclareMathOperator{\inv}{Inv}

\DeclareMathOperator{\iso}{Iso}

\DeclareMathOperator{\orb}{Orb}

\DeclareMathOperator{\poly}{poly}

\DeclareMathOperator{\Sp}{Sp}

\DeclareMathOperator{\reg}{Reg}

\DeclareMathOperator{\rk}{rk}

\DeclareMathOperator{\soc}{Soc}

\DeclareMathOperator{\sym}{Sym}
\DeclareMathOperator{\WL}{WL}

\def\eprf{\hfill$\square$}
					
\def\mmod#1#2#3{#1=#2\ (\text{\rm mod}\hspace{2pt}#3)}
\def\qaq{\quad\text{and}\quad}
\def\qoq{\quad\text{or}\quad}
\def\ov{\overline}
\newcommand{\grp}[1]{\langle {#1}\rangle}

\begin{document}

\title[Testing isomorphism of central Cayley graphs]{Testing isomorphism of central Cayley graphs over almost simple groups in polynomial time}
\author{Ilia Ponomarenko}
\address{St.Petersburg Department of the Steklov Mathematical
	Institute, St.Petersburg, Russia}
\email{inp@pdmi.ras.ru}
\author{Andrey Vasil'ev}
\address{Sobolev Institute of Mathematics, Novosibirsk, Russia\vspace{-7pt}}
\address{Novosibirsk State University, Novosibirsk, Russia}
\email{vasand@math.nsc.ru}
\dedicatory{In memory of Sergei Evdokimov}

\begin{abstract}
A Cayley graph over a group $G$ is said to be central if its connection set is a normal subset of~$G$. It is proved that for any two central Cayley graphs over explicitly given almost simple groups of order $n$, the set of all isomorphisms from the first graph onto the second can be found in time $\poly(n)$. 
\end{abstract}

\maketitle

\section{Introduction}

In the present paper, we are interested in a special case of the following  restriction of the Graph Isomorphism Problem to the class of Cayley graphs.\medskip

{\bf Cayley Graph Isomorphism Problem.} {\it For two explicitly given finite groups $G$ and $G'$ and two sets $X\subset G$ and $X'\subset G'$, construct the set $\iso(\Gamma,\Gamma')$, where $\Gamma=\cay(G,X)$ and $\Gamma'=\cay(G,X')$.}\medskip

Here the input consists of the multiplication tables of $G$ and $G'$ and  the sets $X$ and $X'$, whereas the output is either empty or given by an element of the set $\iso(\Gamma,\Gamma')$ and a generating set of the group $\aut(\Gamma)$ (of size polynomial in the order~$n$ of the group~$G$). Obviously, the Luks algorithm~\cite{L82} solves the Cayley Graph Isomorphism Problem in polynomial time for every group~$G$, whenever the set $X$ is of constant sizes. If $G$ is cyclic, then the problem with no restriction for~$X$ is also solvable in polynomial time (see~\cite{EP03} and~\cite{M01}). It should be noted that if $G=G'$ and $G$ is a so-called CI-group, then an obvious algorithm solves the Cayley Graph Isomorphism Problem in time polynomial in $|\aut(G)|$ (for details, see~\cite{L02}).\medskip

The aforementioned special case is formed by the two following conditions imposed on the input graphs and  groups. First, we assume that $\Gamma$ is a {\it central} Cayley graph over~$G$, which means that $X$ is a normal subset of $G$, i.e., $X^g=X$ for every $g\in G$. Second, the group $G$ is assumed to be {\it almost simple}, i.e., the socle of $G$ is a nonabelian simple group. The same conditions are imposed on the graph $\Gamma'$ and group~$G'$. Even in this rather restrictive case, the problem is still nontrivial; at least the number of possible input graphs is exponential in~$n$.\medskip

{\bf Example.} Let $G=\sym(m)$ be a symmetric group of degree~$m$. Then the number $N(m)$ of central colored Cayley graphs over $G$ is equal to $2^{p(m)}$, where $p(m)$ is the number of all partitions of the set~$\{1,\ldots m\}$. Since $p(m)$ is approximately equal to $2^{\sqrt{m}}$, the number $N(m)$ is exponential in $n=m!$.\medskip

By technical reasons, it is more convenient to deal with {\it colored}  Cayley graphs. Such a graph is given by a partition $\fP$ of the group~$G$ into $k\ge 2$ classes $X_0,\ldots,X_{k-1}$ with $X_0=\{1\}$, and can be thought as arc-colored complete graph with vertex set $G$ and the $i$th color class of arcs coinciding with the arc set of the Cayley graph $\cay(G,X_i)$, $i=0,\ldots,k-1$. We say that $\fP$ is the {\it Cayley partition} of this graph and denote the latter by $\cay(G,\fP)$. In what follows, all Cayley graphs are assumed to be colored: the graph $\cay(G,X)$ is treated as $\cay(G,\fP)$ for $k=3$ and $X_1=X$.

\thrml{140516a}
For any two  central Cayley graphs $\Gamma$ and $\Gamma'$ over explicitly given almost simple groups $G$ and $G'$ of order $n$, the set $\iso(\Gamma,\Gamma')$ can be found in  time $\poly(n)$.
\ethrm

\crllrl{310516y}
The automorphism group of a central Cayley graph over an explicitly given
almost simple group of order $n$ can be found in  time $\poly(n)$.
\ecrllr

The proof of Theorem~\ref{140516a} is a mix of combinatorial and permutation group techniques. Section~\ref{061216a} provides a relevant background for the combinatorial part including coherent configurations and Cayley schemes. In Section~\ref{061216b}, we use a classification of regular almost simple subgroups of primitive groups~\cite{LPS} to prove (Lemma~\ref{081116a}) that except for one special case, if $K\le\sym(G)$ is a $2$-closed primitive group containing regular almost simple subgroup, then
\qtnl{071216a}
K=\sym(G)\qoq K\le D(2,G),
\eqtn
where $D(2,G)$ is the subgroup of $\sym(G)$ generated by the holomorph of~$G$ and the permutation $\sigma:g\mapsto g^{-1}$, $g\in G$. We extend this result to non-primitive groups in Sections~\ref{061216d} and~\ref{251116a} by showing that in this case, either formula~\eqref{071216a} holds or $K$ is a nontrivial generalized wreath product (Theorems~\ref{240416b} and~\ref{240916a}).  We apply this fact in Section~\ref{271116a} to the automorphism group $K$ of a central Cayley graph $\Gamma$ over the group $G$ to establish that
\qtnl{071216b}
K^S=\sym(S)\qoq K^S\le D(2,S),
\eqtn
where $S=\soc(G)$ is the socle of $G$ and  $K^S$ is the restriction to $S$ of the setwise stabilizer of $S$ in $K$. Note that if $G$ is a symmetric group of degree at least~$5$, then the group $D(2,G)$ is isomorphic to the group $G\wr\sym(2)$. Thus as a  byproduct of~\eqref{071216b},  we obtain the following  generalization of the result~\cite[Theorem~1.1]{G16} on the automorphism group of the Cayley graph $\cay(G,X)$, where $G$ is a symmetric group and $X$ is the set of its transpositions.

\thrml{061216e}
Let $G$ be a symmetric group of degree at least~$5$, $X$ a proper normal subset of $G\setminus\soc(G)$, and $\Gamma=\cay(G,X)$. Then $\aut(\Gamma)=D(2,G)$.
\ethrm

In Sections~\ref{061216f} and~\ref{061216g}, we develop algorithmic tools to find the above structure of the group $K$ with the help of the related Cayley scheme and the group~$G$. The main algorithm providing the proof of Theorem~\ref{140516a} is given in Section~\ref{140516w}.\medskip


{\bf Notation.}

The diagonal of the Cartesian product $\Omega\times\Omega$ is denoted by~$1_\Omega$.

For $s\subseteq\Omega\times\Omega$, set $s^*=\{(\beta,\alpha):\ (\alpha,\beta)\in s\}$.

For $\alpha\in\Omega$ and $s\subseteq\Omega\times\Omega$, set $\alpha s=\{\beta\in\Omega:\ (\alpha,\beta)\in s\}$.

For $\Delta\subseteq\Omega$  and $s\subseteq\Omega\times\Omega$, set $s_\Delta=s\cap(\Delta\times\Delta)$.

For a partition $\fE$ of a set $\Omega$ and $s\subseteq\Omega\times\Omega$, set $s_\fE$ to be the relation consisting of the pairs $(\Delta,\Delta')\in\fE\times\fE$ such that $s$ meets $\Delta\times\Delta'$.

For a set $S$ of binary relations, put $S^\cup$ to be the set of all unions of relations from~$S$.

The symmetric and alternating groups on $\Omega$ are denoted by $\sym(\Omega)$ and $\alt(\Omega)$, respectively.

For $f\in\sym(\Omega)$  and $s\subseteq\Omega\times\Omega$, set $s^f=\{(\alpha^f,\beta^f):\ (\alpha,\beta)\in s\}$.

For  a group $G$ and its subgroup $L$, set $L_{left}$ and $L_{right}$ to be the subgroups of $\sym(G)$ induced by left and right multiplications of $L$, respectively, and $L^*=L_{left}L_{right}$.

For a group $G$, set $D(2,G)$ to be the subgroup of $\sym(G)$ generated by the group $\hol(G)=N_{\sym(G)}(G_{right})$ and the permutation $\sigma:g\mapsto g^{-1}$.

For a group $K\le\sym(G)$ and a set $H\subseteq G$, the restriction to $H$ of the setwise stabilizer of $H$ in~$K$ is denoted by $K^H$.

For an imprimitivity system $\fL$ of a transitive group $K$, set $K_\fL$ and $K^\fL$ to be, respectively, the intersection of all $K^H$ with $H\in\fL$ and the permutation group induced by the action of~$K$ on $\fL$.

For a group $G$ and a permutation group $K$, set $\reg(K,G)$ to be the set of all regular subgroups of $K$ that are isomorphic to $G$.

\section{Coherent configurations and Cayley schemes}\label{061216a}
This section contains well-known basic facts on coherent configurations.
All of them can be found in~\cite{EP09} and papers cited there.

\sbsnt{Main definitions.}\label{230415b}
Let $\Omega$ be a finite set and $S$ a partition of~$\Omega\times\Omega$. The pair $\cX=(\Omega,S)$ is called a {\it coherent configuration} on $\Omega$ if the following
conditions hold:
\nmrt
\tm{C1} $1_\Omega\in S^\cup$,
\tm{C2} $S^*=S$,
\tm{C3} given $r,s,t\in S$, the number $c_{rs}^t=|\alpha r\cap\beta s^*|$
does not depend on the choice of the pair $(\alpha,\beta)\in t$.
\enmrt
The elements of $\Omega$ and $S$, and the numbers $c_{rs}^t$ are called the {\it points} and {\it basis relations}, and the {\it intersection numbers} of~$\cX$, respectively. The numbers $|\Omega|$
and $|S|$ are called the {\it degree} and the {\it rank} of~$\cX$. The coherent configuration~$\cX$ is said to be {\it homogeneous} if $1_\Omega\in S$.\medskip

Denote by $\Phi=\Phi(\cX)$ the set of $\Lambda\subseteq\Omega$ such that
$1_\Lambda\in S$. The elements of $\Phi$ are called the {\it fibers}
of~$\cX$. In view of condition (C1), the set $\Omega$ is the disjoint
union of all of them. Moreover, for each $s\in S$, there exist uniquely
determined fibers~$\Lambda$ and~$\Delta$ such that
$s\subseteq\Lambda\times\Delta$. Note that the coherent configuration $\cX$ is homogeneous if and only if $|\Phi|=1$.\medskip

Let $e\in S^\cup$ be an equivalence relation and $\fE$ the set of its classes. Given $\Delta\in\fE$ denote by   $S_\Delta$ the set of all nonempty relations $s_\Delta$ with $s\in S$. Then the pair
$$
\cX_\Delta=(\Delta,S_\Delta)
$$
is a coherent configuration called the {\it restriction} of~$\cX$ to~$\Delta$. This enables to define the restriction of~$\cX$ to a set $\Lambda\in\Phi^\cup$: the corresponding equivalence relation is equal to the union of $\Delta\times\Delta$, where $\Delta$ runs over the fibers contained in~$\Lambda$. Another coherent configuration associated with~$e$ is obtained as follows. Denote by $S_\fE$  the set of all nonempty relations $s_\fE$, $s\in S$. Then
$$
\cX_{\fE}=(\fE,S_{\fE})
$$
is a coherent configuration called the {\it quotient} of $\cX$ modulo~$e$.

\sbsnt{Combinatorial and algebraic isomorphisms.}\label{190117a}
A bijection $f:\Omega\to\Omega'$ is called the (combinatorial) {\it isomorphism} from $\cX$ onto a coherent configuration $\cX'=(\Omega',S')$ if the set $S'$ contains the relation $s^f$ for each $s\in S$. The set of all isomorphisms $f$ is denoted by $\iso(\cX,\cX')$. The group of all isomorphisms of $\cX$ to itself
contains a normal subgroup
$$
\aut(\cX)=\{f\in\sym(\Omega):\ s^f=s,\ s\in S\}
$$
called the {\it automorphism group} of~$\cX$. Conversely, let $G\le\sym(\Omega)$ be a permutation group, and let $S$  be the set of orbits of the component-wise action of $G$
on~$\Omega\times\Omega$. Then the pair $\cX=(\Omega,S)$ is a coherent configuration; we say that $\cX$ {\it is  associated} with~$G$ and denote it by $\inv(G)$.\medskip

According to Wielandt \cite{Wie1969}, a permutation group $G$ on~$\Omega$ is said to be {\it $2$-closed} if it is equal to its {\it $2$-closure}
$$
G^{(2)}=\aut(\inv(G)),
$$
or, equivalently, if $G$ is an automorphism group of a family of binary relations on~$\Omega$ (such a family can always be chosen as the set of basis relations of a coherent configuration on~$\Omega$). If $G$ is $2$-closed and $\fL$ is an imprimitivity system of~$G$, then the group $G_\fL$ is $2$-closed. However, the group $G^\fL$ is not always $2$-closed.\medskip

A bijection $\varphi:S\to S',\ r\mapsto r'$ is called an {\it algebraic isomorphism} from~$\cX$ onto~$\cX'$ if
\qtnl{f041103p1}
c_{rs}^t=c_{r's'}^{t'},\qquad r,s,t\in S.
\eqtn
In this case, $\cX$ and $\cX'$ are said to be {\it algebraically isomorphic}. Each isomorphism~$f$ from~$\cX$ onto~$\cX'$ induces an algebraic isomorphism $\varphi_f:r\mapsto r^f$  between these configurations. The set of all isomorphisms inducing the algebraic isomorphism~$\varphi$ is denoted by $\iso(\cX,\cX',\varphi)$. In particular,
\qtnl{190316b}
\iso(\cX,\cX,\id_S)=\aut(\cX),
\eqtn
where $\id_S$ is the identity mapping on $S$.\medskip

An algebraic isomorphism $\varphi$ induces a bijection from $S^\cup$ onto $(S')^\cup$: the union $r\cup s\cup\cdots$ of basis relations of $\cX$ is taken to $r'\cup s'\cup\cdots$. This bijection is also denoted by $\varphi$. It preserves the equivalence relations $e\in S^\cup$;  moreover, the equivalence relations $e$ and $\varphi(e)$ have the same number of classes as well as the same multiset of their sizes. In this case, if $\fE$ and $\fE'$ are the sets of classes of $e$ and $\varphi(e)$, respectively, and $\Delta\in\fE$, then $\varphi$ induces the algebraic isomorphisms
$$
\varphi_\Delta:\cX^{}_{\Delta^{}}\to\cX'_{\Delta'},\ s\mapsto \varphi(s)_{\Delta'} \qaq\varphi_{\fE}:\cX^{}_{\fE^{}}\to\cX'_{\fE'},\ s_{\fE^{}}\mapsto \varphi(s)_{\fE'}
$$
for a suitable $\Delta'\in\fE'$.

\sbsnt{Direct sum and wreath product.}\label{190117b}
Let $\cX=(\Omega,S)$ and $\cX'=(\Omega',S')$ be coherent configurations. Denote by $\Omega\sqcup\Omega'$ the disjoint union of~$\Omega$ and~$\Omega'$, and by $S\boxplus S'$ the union of the set $S\sqcup S'$ and the set of all relations $\Delta\times\Delta'$ and $\Delta'\times\Delta$ with $\Delta\in\Phi(\cX)$ and $\Delta'\in\Phi(\cX')$. Then the pair
$$
\cX\boxplus\cX'=(\Omega\sqcup\Omega',S\boxplus S')
$$
is a coherent configuration called the {\it direct sum} of~$\cX$ and~$\cX'$. The automorphism group of this configuration equals the direct product $\aut(\cX)\times\aut(\cX')$ acting on the set $\Omega\sqcup\Omega'$. Furthermore, if $\varphi$ is an algebraic isomorphism from $\cX\boxplus\cX'$ to another coherent configuration, then the latter is also the direct sum $\cY\boxplus\cY'$ and~$\varphi$ induces algebraic isomorphisms $\cX\to\cY$ and $\cX'\to\cY'$ coinciding with the restrictions of $\varphi$ on $\Omega$ and $\Omega'$, respectively.\medskip

Let $\cX$ be a homogeneous coherent configuration, $e\in S^\cup$ an equivalence relation, and $\fE$ the set of classes of~$e$. We say that $\cX$ is the {\it wreath product} with respect to~$e$ if for each $s\in S$ such that $s\not\subseteq e$,
$$
s=\bigcup_{(\Lambda,\Delta)\in s_\fE}\Lambda\times\Delta.
$$
In what follows, we always assume that the classes of $e$ can be identified with the help of a family of the isomorphisms $f_{\Lambda,\Delta}:\cX_\Lambda\to\cX_\Delta$, $\Lambda,\Delta\in\fE$, such that
\qtnl{170117a}
(s_{\Lambda})^{f_{\Lambda,\Delta}}=s_{\Delta},\qquad s\in S.
\eqtn
In this case, $\cX$ is isomorphic to the usual wreath product $\cX_\Delta\wr\cX_\fE$ for all $\Delta\in\fE$ (see~\cite[p.45]{W76}). The automorphism group of this coherent configuration is permutation  isomorphic to the wreath product $\aut(\cX_\Delta)\wr\aut(\cX_\fE)$ in imprimitive action.\medskip

Furthermore, if $\varphi$ is an algebraic isomorphism from the wreath product $\cX$ with respect to~$e$ to another coherent configuration, then the latter is also the wreath product $\cX'$ with respect to $e'=\varphi(e)$ and~$\varphi$ induces algebraic isomorphisms $\cX^{}_{\Delta^{}}\to\cX'_{\Delta'}$ and $\cX^{}_{\fE^{}}\to\cX'_{\fE'}$ coinciding with the restrictions of $\varphi$ on $\Delta$ and $\fE$, respectively, where $\fE'$ is the set of classes of~$e'$ and $\Delta'\in\fE'$.

\sbsnt{Cayley schemes.}
A coherent configuration $\cX=(\Omega,S)$ is called the {\it Cayley scheme} over a group~$G$ if
$$
\Omega=G\qaq G_{right}\le\aut(\cX).
$$
In this case, $\cX$ is homogeneous and each basis relation $s$ is the set of arcs of the Cayley graph $\cay(G,X)$, where $X$ is the neighborhood of the identity of $G$ in the relation~$s$. In particular, $\cX$ can be treared as a color graph $\cay(G,\fP)$, where the classes of the Cayley partition~$\fP$ are the neighborhoods of the identity of $G$ in the basis relations of~$\cX$.\medskip

The class of Cayley schemes is closed with respect to taking restrictions and quotients. Namely, if $\cX$ is a Cayley scheme over a group $G$ and $e\in S^\cup$ is an equivalence relation, then the class $H$ of $e$ containing the identity of~$G$ is a subgroup of~$G$. Moreover, the set $\fE$ of classes of~$e$ coincides with the right $H$-cosets of~$G$. It follows that
$$
\cX_H=\cay(H,\fE_H)\qaq \cX_\fE=\cay(G/H,\fE_{G/H}),
$$
with
$$
\fE_H=\{X\in\fE:\ X\subset H\}\qaq\fE_{G/H}=\{\pi(X):\ X\in\fE\},
$$
where in the latter case, $H$ is a normal subgroup of~$G$ and $\pi:G\to G/H$ is the canonical epimorphism.\medskip

Assume that the Cayley scheme $\cX$ is the {\it wreath product} with respect to an equivalence relation~$e\in S^\cup$. Then for any two  classes $\Lambda,\Delta\in\fE$, there exists a permutation $f\in G_{right}$ taking $\Lambda$ to $\Delta$; set $f_{\Lambda,\Delta}$ to be the restriction of $f$ to~$\Lambda$. Since all these $f$ are automorphisms of $\cX$, the family $\{f_{\Lambda,\Delta}\}_{\Lambda,\Delta\in\fE}$ satisfies conditions~\eqref{170117a}.\medskip

The Cayley scheme $\cX$ is said to be {\it central} if $G_{left}\le\aut(\cX)$ which is by definition of a Cayley scheme is  equivalent to $G^*\le\aut(\cX)$. One can see that $\cX$ is central if and only if the colored Cayley graph associated with $\cX$ is central.

\sbsnt{Partial order and the WL-algorithm.}
There is a natural partial order\, $\le$\, on the set of all coherent configurations on the same set~$\Omega$. Namely, given two coherent configurations $\cX=(\Omega,S)$ and
$\cX'=(\Omega,S')$, we set
$$
\cX\le\cX'\ \Leftrightarrow\ S^\cup\subseteq (S')^\cup.
$$
The minimal and maximal elements with respect to this order are, respectively, the {\it trivial} and {\it discrete} coherent configurations. The first one is a unique coherent configuration $\cT_\Omega$ with at most two basis relations: $1_\Omega$ and its complement to $\Omega\times\Omega$ (if $\Omega$ consists of at least two points). Every basis relation of the discrete configuration is a singleton. With respect to this order, the direct sum $\cX\boxplus\cX'$ is the minimal coherent configuration on $\Omega\sqcup\Omega'$, the restrictions of which to $\Omega$ and $\Omega'$ are equal to $\cX$ and $\cX'$, respectively.\medskip

One can prove that given a set $T\subseteq 2^{\Omega\times\Omega}$, there exists a unique minimal coherent configuration $\cX'$ such that every relation of~$T$ is the union of some basis relations of~$\cX'$. This coherent configuration is called the {\it coherent closure} of $T$ and can be constructed by the well-known Weisfeiler-Leman algorithm (WL-algorithm) \cite[Section~B]{W76} in time polynomial in sizes of~$T$ and~$\Omega$. To stress this fact,  the coherent closure of $T$ is denoted by $\WL(T)$. For a color graph $\Gamma$ with the set $S$ of color classes,  we set
$$
\WL(\Gamma,T)=\WL(S\cup T)
$$
and write $\WL(\Gamma)$ instead of $\WL(\Gamma,\varnothing)$.
It is important to note that the automorphism group of the coherent configuration $\WL(\Gamma,T)$ is equal to the subgroup of $\aut(\Gamma)$ leaving each relation of~$T$ fixed (as a set). This implies that if $\Gamma$ is a Cayley graph over $G$, then the coherent configuration $\WL(\Gamma)$ is a Cayley scheme over $G$. Since any coherent configuration $\cX$ can be considered as a color graph, we extend our notation to write $\WL(\cX,T)$. Concerning the following statement, we refer to~\cite[Theorem~2.4]{P13}.

\thrml{131116t}
Let $S$ and $S'$ be $m$-sets of binary relations on an $n$-element set. Then given a bijection $\psi:S\to\cS'$, one can check in time $mn^{O(1)}$ whether or not there exists an algebraic isomorphism $\varphi:\WL(S)\to\WL(S')$ such that $\varphi|_S=\psi$. Moreover,
if $\varphi$ does exist, then it can be found within the same time.
\ethrm

\section{Almost simple groups}\label{061216b}

In this section, we collect several known facts on finite almost simple groups and deduce some auxiliary results to be used throughout the paper.

\lmml{310516a}
Let $G$ be an almost simple group of order~$n$. Then
\nmrt
\tm{i} $|G/L|\le \log n$, where $L=\soc(G)$,
\tm{ii} $|\reg(K,G)|=n^{O(1)}$ for every $K\le\sym(G)$ containing $G$ as a normal regular subgroup.
\enmrt
\elmm
\proof From the description of the automorphism groups of simple groups (see, e.g., \cite[Introduction]{CCNPW}), it follows that $|\aut(L)/L|\le \log|L|$. Therefore, statement~(i) is a consequence of the inclusions $L\le G\le\aut(G)\le\aut(L)$. The inclusions also imply that $|\aut(G)|\le n\log n$. Since the centralizer of $G$ in $\sym(G)$ is of order $n$ \cite[Exercise 4.5']{Wie1964}, we have
$$
|K|\le |N_{\sym(G)}(G)|\le |C_{\sym(G)}(G)|\cdot|\aut(G)|\le n^2\log n.
$$
On the other hand, the group $G$ is $3$-generated \cite{DL95}. Thus, the number of regular subgroups of $K$ isomorphic to $G$ is at most $\binom{K}{3}$ and statement (ii) follows.\eprf\medskip

In the following statement, we use the classification of regular almost simple subgroups of a primitive group~\cite[Theorem~1.4]{LPS}.

\lmml{081116a}
Let $G$ be an almost simple group and $G^*\le K\le\sym(G)$. Suppose that $K$ is primitive. Then one of the following holds:
\nmrt
\tm{i} $K\ge\alt(G)$,
\tm{ii} $G=\soc(G)$ and $K\le D(2,G)$,
\tm{iii} $G=\sym(5)$ and $K=Sp_8(2)$ is $2$-transitive.
\enmrt

\elmm
\proof Without loss of generality, we assume that neither $K\ge\alt(G)$, nor $G=\soc(G)$ and $K\le D(2,G)$. Then by aforementioned classification, exactly one of the following pairs $(G,\soc(K))$ occurs:
\nmrt
\tm{a} $(\alt(p^2-2),\alt(p^2+1))$, where $\mmod{p}{3}{4}$ is prime,
\tm{b} $(\sym(p-2),\alt(p))$ or $(\sym(p-2),\alt(p+1))$, where $p\ge 7$ is
prime,
\tm{c} the twelve pairs in the table below.
\enmrt
\medskip
\centerline{\begin{tabular}{|c|c|c|}
	\hline
	\rule[-1ex]{0pt}{2.5ex} $No.$ & $G$ & $\soc(K)$ \\
	\hline
	\rule[-1ex]{0pt}{2.5ex} $1$   & $\alt(5)$ & $L_2(59)$ \\
	\hline
	\rule[-1ex]{0pt}{2.5ex} $2$   & $\alt(7)$ & $\alt(11)$ \\
	\hline
	\rule[-1ex]{0pt}{2.5ex} $3$   & $\alt(7)$ & $\alt(12)$  \\
	\hline
	\rule[-1ex]{0pt}{2.5ex} $4$   & $\sym(5)$ & $\alt(9)$ \\
	\hline
	\rule[-1ex]{0pt}{2.5ex} $5$   & $\sym(5)$ & $\Sp_4(4)$ \\
	\hline
	\rule[-1ex]{0pt}{2.5ex} $6$   & $\sym(5)$ & $\Sp_6(2)$ \\
	\hline
\end{tabular} \qquad
\begin{tabular}{|c|c|c|}
	\hline
	\rule[-1ex]{0pt}{2.5ex} $No.$ & $G$ & $\soc(K)$ \\
	\hline
	\rule[-1ex]{0pt}{2.5ex} $7$   & $\sym(5)$ & $\Omega_8^+(2)$ \\
	\hline
	\rule[-1ex]{0pt}{2.5ex} $8$   & $\sym(5)$ & $\Sp_8(2)$ \\
	\hline
	\rule[-1ex]{0pt}{2.5ex} $9$   & $L_2(16).4$ & $\Sp_6(4)$  \\
	\hline
	\rule[-1ex]{0pt}{2.5ex} $10$  & $L_2(16).4$ & $\Omega_8^+(4)$ \\
	\hline
	\rule[-1ex]{0pt}{2.5ex} $11$  & $L_3(4).2$ & $M_{23}$ \\
	\hline
	\rule[-1ex]{0pt}{2.5ex} $12$  & $L_3(4).2$ & $M_{24}$ \\
	\hline
\end{tabular}}\medskip\medskip

\noindent The assumption $G^*\le K$, in particular, implies that $|G|^2\le |K|$. By straightforward check this excludes cases (a) with $p>3$, (b) with $p>7$, and cases 11 and 12 from the table. Similarly, the remaining cases in (a) and (b) as well as cases 1--4, and 9 are impossible because $|G|^2$ must divide $|K|$.\medskip

In cases 5, 6, 7, and 10 from the table, we check the maximal subgroups of $K$ and show that none of them includes the subgroup isomorphic to $G^*$. Indeed, in cases 5 and 6 none of the maximal subgroups contains $\alt(5)\times\alt(5)$ \cite[Tables~8.14, 8.28, 8.29]{BHR}. In case 10, information from \cite[Table~2]{LPS} shows that $K$ is an extension of $\Omega_8^+(4)$ by a field automorphism. This group contains the only (up-to conjugation) maximal subgroup with section isomorphic to $L_2(16)\times L_2(16)$ \cite[Table~8.50]{BHR}, but the order of this subgroup is less than $|G^*|$. In case 7, we make use of \cite[Table~2]{LPS} to see that $K=\Omega_8^+(2)$. Again this group includes up to conjugation the only maximal subgroup $M$ with section isomorphic to $\alt(5)\times\alt(5)$ (see \cite[Table~8.50]{BHR}). However, $|M|=|G^*|$ but $M\simeq D(2,\alt(5))$ is not isomorphic to $G^*\simeq\sym(5)\times\sym(5)$.\medskip

This leaves us with case 8 of the table where $\soc(K)=Sp_8(2)=K$ and we arrive at case (iii) of the conclusion of the lemma.\eprf\medskip

\section{The structure of automorphism groups: the principal section}\label{061216d}

\sbsnt{Preliminaries.}\label{251116t}
Let $G$ be a finite group. The automorphism group of every central Cayley graph over $G$ contains a subgroup $G^*$ (see Notation). In this section, we establish some basic facts on the permutation groups $K$ satisfying the following condition:
\qtnl{240416a}
G^*\le K\le\sym(G),
\eqtn
where $G$ is an almost simple group. We use a concept of the {\it generalized wreath product} of permutation groups introduced and studied in~\cite{EP12}. Namely, a transitive group~$K$ is the generalized wreath product if it has two imprimitivity systems $\fL$ and $\fU$ such that every block of $\fL$ is contained in a block of~$\fU$ and
$$
K_\fL=\prod_{X\in\fU}(K_\fL)^X.
$$
The generalized wreath product is said to be {\it trivial} if either $\fL$ consists of singletons or $\fU=\{\Omega\}$. When $\fL=\fU$, the group $K$ is permutation isomorphic to the wreath product $K^X\wr K^\fL$ in imprimitive action, where $X\in\fU$. 

\thrml{240416b}
Let $G$ be an almost simple group, and let $K$ be a $2$-closed group
satisfying condition~\eqref{240416a}. Then one of the following statements
holds:
\nmrt
\tm{i} $K=\sym(G)$ or $G^*\trianglelefteq K$,
\tm{ii} $K$ is a nontrivial generalized wreath product.
\enmrt
\ethrm

The proof of Theorem~\ref{240416b} is given in the end of Section~\ref{271116a}.\medskip

\sbsnt{The minimal block.}\label{290516a}
Let $K$ satisfy condition~\eqref{240416a} and $X$ a $K$-block containing the identity of~$G$. Since $K$ is a permutation group on $G$ that includes $G_{right}$, the block $X$ is a subgroup of $G$ \cite[Theorem~24.12]{Wie1964}. Taking into account that $G_{left}$ also lies in~$K$, we conclude that $X$ is normal. Denote by~$L$ the intersection of all non-singleton $K$-blocks containing the identity of the group~$G$. Then $L$ is a $K$-block and we call it the {\it minimal block} of $K$.

\lmml{081116d}
Let $K$ satisfy condition~\eqref{240416a} and $L$ the minimal block of $K$. Then $L$ is normal subgroup of $G$ including~$\soc(G)$. In particular, $L$ is an almost simple group such that  $\soc(L)=\soc(G)$.
\elmm
\proof According to the above remark, every $K$-block containing the identity of~$G$ is a normal subgroup of~$G$. If the block is not a singleton, then this normal subgroup is nontrivial and hence contains~$\soc(G)$, because the group $G$ is almost simple. Thus, the minimal block $L$ being the intersection of nontrivial normal subgroups of $G$ is a normal subgroup and contains $\soc(G)$.\eprf\medskip

Denote by $\fL$ the imprimitivity system containing~$L$. Obviously,
$$
\orb(L^*,G)=\fL.
$$
Recall that according to the definition, $K_\fL$ is a normal subgroup of~$K$ leaving each block of~$\fL$ fixed (as a set), and given $X\in\fL$, the group $(L^*)^X\le\sym(X)$ is induced by right multiplications of~$L$.

\lmml{071116a}
For any $X\in\fL$, the group $(K_\fL)^X$ is primitive and contains $(L^*)^X$.
\elmm
\proof The normality of $L$ in $G$ implies that the orbits of the action of $L^*$ on $G$ coincide with the $L$-cosets. It follows that the permutation group induced by this action is contained in~$K_\fL$. This proves the second statement. To prove the first statement, in view of the transitivity of $K$, we may assume that $X=L$.\medskip

Assume on the contrary that the group $K'=(K_\fL)^L$ is not primitive.  Then there exists a minimal non-singleton $K'$-block $L'<L$. Taking into account that $K'\trianglelefteq K^L$, we conclude that $L''=(L')^k$ is a also $K'$-block for every $k\in K$ \cite[Proposition~6.2]{Wie1964}. The imprimitivity system $\fL'$ of the group~$K'$ that contains~$L'$ coincides with the imprimitivity system containing~$L''$, for otherwise by the minimality of~$L'$ and Lemma~\ref{081116d} applied for the group $K'$, one can choose the block $L''$ so that
$$
1=L'\cap L''\supseteq \soc(L)=\soc(G)\supsetneq 1,
$$
a contradiction. Thus, $\fL'$ is an imprimitivity system of the group $K^L$. Consequently, $L'$ is a non-singleton $K$-block strictly contained in $L$, which is impossible by the definition of~$L$.\eprf

\sbsnt{The wreath decomposition of~$K_\fL$.}\label{171116v}
For every two sets $X,X'\in\fL$, we write $X\sim X'$ if the restriction epimorphisms
\qtnl{071116y}
(K_\fL)^{X\cup X'}\to(K_\fL)^{X^{}}\qaq (K_\fL)^{X\cup X'}\to(K_\fL)^{X'}
\eqtn
are isomorphisms. In particular, the groups $(K_\fL)^{X^{}}$ and $(K_\fL)^{X'}$ are isomorphic. It is easily seen that $\sim$ is an equivalence relation on~$\fL$. This relation is $K$-invariant, because $K_\fL$ is a normal subgroup of~$K$. Denote by $U$ the union of $L$-cosets belonging to the class of $\sim$ that contains $L$. Then $U$ is obviously a $K$-block and hence is a normal subgroup of $G$. Thus,
$$
L\trianglelefteq U\trianglelefteq G.
$$
The imprimitivity system of the group~$K$ that contains the block~$U$ is denoted by~$\fU$. We say that $U/L$ is the {\it principal} $K$-section of~$G$, and $\fU$ and $\fL$ are the {\it associated} partitions.

\thrml{200416b}
Let $G$ be an almost simple group, $K$ a $2$-closed group satisfying condition~\eqref{240416a}, and $\fU$, $\fL$ are the partitions associated with the principal $K$-section. Then the coherent configuration $\inv(K_\fL)$ is the direct sum of the coherent configurations $\inv(K_\fL)_Y$, where $Y\in\fU$. In particular,
\qtnl{240516a}
K_\fL=\prod_{Y\in\fU}(K_\fL)^Y,
\eqtn
i.e., $K$ is generalized wreath product.
\ethrm
\proof The subgroup $K_\fL$ of $2$-closed group $K$ is $2$-closed too (see Subsection~\ref{190117a}). Therefore, equality \eqref{240516a} follows from the first statement of the theorem, because the automorphism group of the direct sum equals the direct product of the summands (see  Subsection~\ref{190117b}). To prove the first statement, it suffices to verify that given $X,X'\in\orb(K_\fL,G)$,
\qtnl{020616a}
X\not\sim X'\quad\Rightarrow\quad X\times X'\in\orb(K_\fL,X\times X').
\eqtn
To this end, we note that the group $(K_\fL)^{X\cup X'}$ is the subdirect product of the transitive constituents $M=(K_\fL)^{X^{}}$ and $M'=(K_\fL)^{X'}$. Therefore, there exist uniquely determined normal subgroups $H$ and $H'$ of $M$ and $M'$, respectively, and a group isomorphism $\varphi:M/H\to M'/H'$ such that
\qtnl{020616y}
(K_\fL)^{X\cup X'}=\{(k,k')\in M\times M':\ \varphi(k)=k'\}.
\eqtn
Now if $X\not\sim X'$, then at least one of the epimorphisms~\eqref{071116y} is not an isomorphism. Therefore, one of the groups $H$ and $H'$, say $H$, is nontrivial. It follows that $H$ being a normal subgroup of the primitive group $M$ (Lemma~\ref{071116a}) acts transitively on $X$. By~\eqref{020616y}, this implies that $(K_\fL)^{X\cup X'}$ contains the subgroup $H\times 1$. Thus,
$$
(x,x')^{H\times 1}=X\times \{x'\}
$$
for all $x'\in X'$ and hence the group $(K_\fL)^{X\cup X'}$ is transitive on the set $X\times X'$.\eprf

\section{The normalizer of $\soc(G)^*$ in $\sym(G)$}\label{251116a}

The goal of this section is to prove the following theorem that shows (as we will see) that the case $U=G$ is very similar to the case where the group $K$ is primitive.

\thrml{240916a}
Let $G$ be an almost simple group, $S=\soc(G)$, and $N=N_{\sym(G)}(S^*)$. Then
\qtnl{231116a}
N=D(2,G).
\eqtn
\ethrm

Clearly, $N$ is a proper subgroup of $\sym(G)$ and condition~\eqref{240416a} is satisfied for $K=N$.
In particular,  the minimal block of $N$ coincides with $S$.

\lmml{231116d}
In the above notation, $C_{\sym(G)}(S^*)=1$.
\elmm
\proof Set $C=C_{\sym(G)}(S^*)$. Then obviously, $C=C_N(S^*)$. Since $N$ is transitive and $C$  is normal in $N$, the orbits of $C$ form an imprimitivity system of~$N$. Denote by $X$ the block of this system that contains the identity of $G$. We may assume that the block $X$ is not a singleton, for otherwise $C=1$ and we are done. Then $S\subseteq X$, because $S$ is the minimal $N$-block. Since $C$ is transitive on $X$ and $S$ is a block of $C^X$, there exists $c\in C$ such that $S^c=S$ and $c^S\ne 1$. However, the latter is impossible, because $c^S$ centralizes the subgroups $S_{left}$ and $S_{right}$.\eprf\medskip

{\bf Proof of Theorem~\ref{240916a}.}  By Lemma~\ref{231116d}, there exists a monomorphism from  $N=N_{\sym(G)}(S^*)$ to $\aut(S^*)$. Since the latter is isomorphic to the wreath product $W=\aut(S)\wr\sym(2)$, this monomorphism induces a monomorphism
$$
\varphi:N\to W.
$$
Clearly, $\varphi$ can be chosen so that the subgroups $S_{left}$ and $S_{right}$ of the group~$N$ go to the subgroup $\Inn(S)\times 1$ of the group $A=\aut(S)\times 1$ and to the subgroup $1\times\Inn(S)$ of the group $B=1\times \aut(S)$, respectively. Set
$$
W_1=\im(\varphi)\qaq W_0=\im(\varphi)\cap(A\times B).
$$
Then since the index of $A\times B$ in $W$ equals~$2$ and $W_1$ is not contained in $A\times B$, we conclude  that
\qtnl{231116c}
|W_1:W_0|=2.
\eqtn
Note that the centralizer of the group $\varphi(S_{right})\le B$ in the group $W_0$ is contained in $A$, because the group $B$ is almost simple with the socle $\varphi(S_{right})$. Since the group $G_{left}\le N$ centralizes~$S_{right}$, this proves the first of the two following inclusions (the second one can be proved in a similar way):
$$
\varphi(G_{left})\le A \qaq \varphi(G_{right})\le B.
$$
The first inclusion implies that $\varphi(G_{left})\le W_0\cap  A$. The reverse inclusion follows from the fact that the centralizer of $G_{left}$ in $\sym(G)$ is equal to $G_{right}$   \cite[Proposition~4.3]{Wie1964}. Thus, we obtain the equalities:
\qtnl{231116b}
\varphi(G_{left})=W_0\cap A\qaq \varphi(G_{right})=W_0\cap B.
\eqtn
This immediately implies that $\varphi(G_{left})$ and $\varphi(G_{right})$  are normal in~$W_0$. This group has trivial center and hence can be identified with a subgroup of the direct product of the groups $A'=\aut(\varphi(G_{left}))$ and $B'=\aut(\varphi(G_{right}))$ (isomorphic to $\aut(G)$). It follows that
$$
\ov{W}_0=W_0/\varphi(G^*)\le\ov{A}'\times\ov{B}',
$$
where $\ov{A}'=A'/\varphi(G_{left})$ and $\ov{B}'=B'/\varphi(G_{right})$. Moreover in view of formulas~\eqref{231116b}, the group $\ov{W}_0$ intersects each of the groups $\ov{A}'$ and $\ov{B}'$ trivially. Therefore,
\qtnl{280317a}
|\ov{W}_0|=|\aut(G)/G|.
\eqtn
Now, the inclusion $\varphi(\aut(G))\le W_0$ and formula~\eqref{280317a} show that $W_0=\varphi(\hol(G))$. Formula~\eqref{231116c} yields that $|N:\hol(G)|=2$. Since the permutation $\sigma:g\mapsto g^{-1}$ lies in $N$, formula~\eqref{231116a} holds.\eprf

\section{Symmetric and normal types of the automorphism group}\label{271116a}

Let $G$ be an almost simple group, $K$ be a $2$-closed group satisfying condition~\eqref{240416a}, and $U/L$ the principal $K$-section of~$G$. Then the group $L$ is almost simple (Lemma~\ref{081116d}) and $K^L$ is primitive (Lemma~\ref{071116a}). We say that
$K$ is of {\it symmetric type} if either $K^L\ge\alt(L)$, or $G=\sym(5)$ and $K=\Sp_8(2)$ is $2$-transitive; if $L=\soc(L)$ and $K^L\le D(2,L)$, the group $K$ is said to be of  {\it normal type}. The following statement is a straightforward consequence of Lemma~\ref{081116a} and the above definitions.

\prpstnl{111116b}
Let $G$ be an almost simple group, and let $K\le\sym(G)$ be a $2$-closed group containing~$G^*$. Then $K$ is of symmetric or normal type.
\eprpstn

Let us study a group of symmetric and normal types in detail. As the following statement shows, any group of symmetric type is, in fact, the wreath product in imprimitive action. In what follows, $\fU$ and $\fL$ are the  partitions associated with the principal $K$-section of~$G$.

\thrml{091116a}
Let $K$ be a group of symmetric type. Then
$$
\fL=\fU\qaq(K_\fL)^L=\sym(L).
$$
In particular, $K$ is permutation isomorphic to the wreath product $\sym(L)\wr K^\fL$ in imprimitive action.
\ethrm
\proof Let $X,X'\in\fL$ and $X\sim X'$. We claim that there exists a bijection  $f:X\to X'$, for which
\qtnl{101116a}
s_f\in\orb(K_\fL,X\times X'),
\eqtn
where $s_f=\{(\alpha,\alpha f): \alpha\in X\}$ is the graph of~$f$. Indeed, consider the group $M=(K_\fL)^{X\cup X'}$. Since $X\sim X'$, the group $M$ acts faithfully on $X$ and $X'$. As the group $K$ is of symmetric type, each of these actions is 2-transitive. However, $M$ has a unique faithful 2-transitive representation of degree $d=|X|$: this is obvious if $M\ge\alt(d)$ and follows from the classification of 2-transitive groups if $M=\Sp_8(2)$ (see, e.g., \cite[Table~7.4]{Cam}). Consequently, among $2d$ point stabilizers $M_\beta$, $\beta\in X\cup X'$, there are exactly $d$ distinct, and also $M_\beta\ne M_\gamma$ whenever $\beta$ and $\gamma$ are distinct points in $X$. Therefore, for every $\alpha\in X$ there is the only $\alpha'\in X'$ such that
$$
M_\alpha=M_{\alpha'}=M_{\alpha,\alpha'}.
$$
Thus, the required bijection $f$ takes $\alpha$ to $\alpha'$.\medskip

To prove that $\fL=\fU$, assume on the contrary that $U$ contains a block $X\in\fL$ other than~$L$. Denote by $f:L\to X$ the bijection defined in the above claim for $X=L$ and $X'=X$. Then by the assumption, the element $f(1)$, where $1$ is the identity of the group $G$, does not belong to $L$. The binary relation~$s_f$ is invariant with respect to the group $K_\fL\ge L^*$  and hence for all $x\in L$,
$$
(x,f(1)x)=(1,f(1))^{x_r}=(x,f(x))=(1,f(1))^{x_l}=(x,xf(1)),
$$
where $x_r$ and $x_l$ are the permutations of $L^*$ induced by the right and left multiplication by~$x$. This implies that the element $f(1)\in X\subset U$ centralizes $L$. However, this is impossible, because $f(1)\ne 1$ and the group $U$ is almost simple.\medskip

Let us prove the second equality. The $2$-closedness of $K$ implies that $K_\fL$ is $2$-closed. Furthermore, in view of the $2$-transitivity of the group $(K_\fL)^X$ its   the $2$-closure equals $\sym(X)$. By Theorem~\ref{200416b} and equality $\fL=\fU$, this implies that
$$
K_\fL=(K_\fL)^{(2)}=(\prod_{X\in\fL}(K_\fL)^X)^{(2)}=\prod_{X\in\fL}((K_\fL)^X)^{(2)}=\prod_{X\in\fL}\sym(X).
$$
Thus, $(K_\fL)^X=\sym(X)$ for all $X\in\fL$ and we are done.\eprf\medskip

Theorem~\ref{091116a} shows that in the case of symmetric type, the group $K^U=\sym(U)$ is the largest possible. In the normal type case, the group $K^U$ is quite small. More exactly, the following statement holds.

\thrml{091116y}
Let $K$ be a group of normal type. Then $K^U\le D(2,U)$.
\ethrm
\proof  By the hypothesis of the theorem, $L=\soc(L)$ and $K^L\le D(2,L)$. By Lemma~\ref{081116d}, the first equality implies that $L=\soc(G)$ and hence
\qtnl{091116f}
L=\soc(U).
\eqtn
The second inclusion implies that $(L^*)^X$ is  a characteristic subgroup of the group $(K_\fL)^X$ for all $X\in\fL$. By the definition of~$U$, this implies that $(L^*)^U$ is a characteristic subgroup of the group $(K_\fL)^U$. However, the latter group is normal in~$K^U$. Thus,
$$
(L^*)^U\trianglelefteq K^U.
$$
It follows that $K^U$ is contained in the normalizer of $(L^*)^U$ in $\sym(U)$. However, this normalizer is contained in $D(2,U)$ by Theorem~\ref{240916a} applied for $G=U$ with taking into account equality~\eqref{091116f}.\eprf\medskip

{\bf Proof of Theorem~\ref{240416b}.}
Let $U/L$ be the principal $K$-section of the group~$G$. By Proposition~\ref{111116b}, the group $K$ is of symmetric or normal type. Suppose first that $U=G$. Then statement~(i) of Theorem~\ref{240416b} holds. Indeed, if $K$ is of symmetric type, then $K=\sym(G)$  (Theorem~\ref{091116a}), whereas if $K$ is of normal type, then $K\le D(2,G)$ (Theorem~\ref{091116y}) and the required statement follows from the fact that $G^*$ is normal in $D(2,G)$. Finally, if $U<G$, then statement~(ii) of Theorem~\ref{240416b} holds by Theorem~\ref{200416b}.\eprf\medskip

{\bf Proof of Theorem~\ref{061216e}.} First, assume that the group $K=\aut(\Gamma)$ is a nontrivial generalized wreath product. Note that $S=\soc(G)$ is a unique proper normal subgroup in~$G$  and $|G/S|=2$. Therefore the generalized wreath product must be a usual one and  the group~$K$ is permutation isomorphic to the wreath product $M\wr C_2$ in imprimitive action for some group $M\le\sym(S)$. It follows that $G\setminus S$ is an orbit of the point stabilizer $K_\alpha$, where $\alpha$ is the identity of~$G$. Since $X\subset G\setminus S$ is a union of some orbits of $K_\alpha$, we conclude that $X=G\setminus S$, a contradiction.\medskip 

The group $K$ is 2-closed as the automorphism group of a graph. The normality of $X$ implies that $K$  satisfies condition~\eqref{240416a} with $G=\sym(m)$ for $m\ge 5$. Finally, $K$ is not a nontrivial generalized wreath product by above, and $K\ne\sym(G)$, because the graph $\Gamma$ is neither complete nor empty. Thus, by Theorem~\ref{240416b}, the group $G^*$ is normal in~$K$ and
$$
G^*\le K\le D(2,G),
$$
hence $K=G^*$ or $K=D(2,G)$, because $|D(2,G):G^*|=2$. However, since $X$ is a normal subset of a symmetric group, we have $X=X^{-1}$, so the graph $\Gamma$ has the automorphism $\sigma:g\mapsto g^{-1}$, $g\in G$. Thus, $K=\grp{G^*,\sigma}=D(2,G)$, as required.\eprf

\section{Finding the principal section in a Cayley scheme}\label{061216f}
\sbsnt{The main resut.}
Let $\cX$ be a central Cayley scheme over an almost simple group $G$. Then the group $K=\aut(\cX)$ is $2$-closed and satisfies condition~\eqref{240416a}. Therefore, $K$ is of symmetric or normal type by Proposition~\ref{111116b}. In this section, we develop an algorithmic technique to determine (with no $K$ in hand) which of these cases occurs for the scheme~$\cX$. The main result here is Theorem~\ref{121116b} below which immediately follows from Corollaries~\ref{111116g} and~\ref{121116y} proved in Subsections~\ref{270516b} and~\ref{121116x}, respectively.

\thrml{121116b}
Given a central Cayley scheme $\cX$ over an almost simple group~$G$ of order $n$, one can determine the type of~$K=\aut(\cX)$ and find the principal $K$-section of $G$ in time $\poly(n)$.
\ethrm

\sbsnt{The case of symmetric type.}\label{270516b}
For a group $H\le G$, denote by $\fH$  and  $\WL(\cX,\fH)$ the partition of $G$ into the right $H$-cosets and the Cayley scheme $\WL(\cX,T)$ with $T=\{1_X:\ X\in\fH\}$, respectively. Recall that $\cT_X$ is the trivial coherent configuration on~$X$. Denote by $\cH_0$ the set of groups $H$ such that $\soc(G)\le H\leq G$ and
\qtnl{111116t}
\WL(\cX,\fH)={{\underset{X\in\fH}{\boxplus}}}\cT_X.
\eqtn

\lmml{111116e}
In the above notation, the following statements hold:
\nmrt
\tm{i} if the group $K$ is of symmetric type, then the set $\cH_0$ is nonempty and the minimal block $L$ of~$K$ is the largest (by inclusion) element of $\cH_0$,
\tm{ii} $K$ is of symmetric type if and only if $\cH_0$ contains $S=\soc(G)$,
\enmrt
\elmm
\proof To prove statement (i), assume that the group $K$ is of symmetric type. Set $\fU$ and $\fL$ to be the partitions associated with the principal $K$-section of~$G$. Then by Theorem~\ref{091116a}, we have $\fL=\fU$ and $\inv(K^X)=\cT_X$ for all $X\in\fL$. By Theorem~\ref{200416b}, this implies that
\qtnl{290317a}
\inv(K_\fL)={{\underset{X\in\fL}{\boxplus}}}\cT_X.
\eqtn
The minimality of the direct sum implies that  $\WL(\cX,\fL)=\inv(K_\fL)$, which proves formula~\eqref{111116t} for $H=L$ and $\fH=\fL$, in particular, $\cH_0$ is nonempty. If $L$ is not the largest element of $\cH_0$, then there exists $H\in\cH_0$ such that $H\setminus L\ne\varnothing$. It follows that $$K\ge\sym(H)\times \id_{G\setminus H},$$ where $\id_{G\setminus H}$ is the identity subgroup of $\sym(G\setminus H)$. Hence there is a permutation $k\in K$ that moves the identity of $G$ to~$H\setminus L$ and leaves all non-identity elements of $L$ fixed. But this is impossible, because $L$ is a $K$-block.\medskip

To prove the necessity for statement~(ii), let $K$ be of symmetric type. Then formula~\eqref{290317a} holds. Therefore, if $\fS$ is the partition of $G$ into the cosets of~$S$, then $\fS$ refines $\fL$ and hence
$$
\WL(\cX,\fS)=\WL(\WL(\cX,\fL),\fS)={{\underset{X\in\fL}{\boxplus}}}\WL(\cT_X,\fS_X)={{\underset{X\in\fS}{\boxplus}}}\cT_X,
$$
where $\fS_X$ is the partition of $X$ induced by $\fS$. Consequently, $S\in\cH_0$. Conversely, assume on the contrary that $K$ is of normal type. Then $K^U\le D(2,U)$ by Theorem~\ref{091116y}. Therefore $K^S\le\hol(S)$. On the other hand, since $S\in\cH_0$, we have $\sym(S)=K^S\le\hol(S)$, a contradiction.\eprf\medskip

From statement~(i) of Lemma~\ref{310516a}, it follows that the number of groups~$H\le G$ containing $\soc(G)$ is at most $\log n$, where $n$ is the order of $G$; in particular, $|\cH_0|\le \log n$. Moreover, for each $H$,  the coherent configuration $\WL(\cX,\fH)$ can be efficiently found by the WL-algorithm and condition~\eqref{111116t}
can be verified by checking at most $|\fH|^2\le n^2$ basis relations. Therefore, the set $\cH_0$ can be found in time $\poly(n)$. By statement~(ii) of Lemma~\ref{111116e}, this is enough to test efficiently whether or not $K$ is of symmetric type, and if it so, then to find the minimal block $L=U$ of the group~$K$ (statement~(i) of the same lemma).

\crllrl{111116g}
Given a central Cayley scheme $\cX$ over an almost simple group $G$ of order $n$, one can test in time $\poly(n)$ whether the group $\aut(\cX)$ is of symmetric type, and (if so) find the principal section of~$\aut(\cX)$ within the same time.
\ecrllr

\sbsnt{The case of normal type.}\label{121116x}
In view of Corollary~\ref{111116g} and Proposition~\ref{111116b}, one can efficiently test whether the automorphism group $K$ of a central Cayley scheme $\cX$ is of normal type. 
Denote by $\cH_1$ the set of all groups $H$ such that $\soc(G)\le H\trianglelefteq G$ and
\qtnl{131116a}
H^*\times \id_{G\setminus H}\le K,
\eqtn
where the left-hand side denotes the subgroup of $\sym(G)$ that leaves each point of $G\setminus H$ fixed and coincides with $H^*$ on $H$.

\lmml{121116a}
Suppose that the group  $K$ is of normal type and $U/L$ is the principal $K$-section of~$G$. Then $L=\soc(G)$, the set $\cH_1$ is nonempty, and $U$ is the smallest element of~$\cH_1$.
\elmm
\proof Lemma~\ref{081116a} yields that $L=\soc(G)$. From Theorem~\ref{200416b}, it is easily follows that $U\in\cH_1$. Assume on the contrary that the group $U$ is not the smallest in $\cH_1$. Then there is a group $V\in\cH_1$ such that $L\le U\cap V<U$. Take a non-identity element $w\in U\cap V$ and denote by $k_U$ (respectively, $k_V$) the permutation on $G$ acting on $U$ (respectively, $V$) by right multiplication by $w$ and acting trivially outside $U$ (respectively, $V$). Then $k_U,k_V\in K$,  the permutation  $$k=k_U^{}k_V^{-1}$$ is not identity on $U$, and $U^k=U\qaq k^{U\cap V}=\id_{U\cap V}.$
Clearly, $k^U$ belongs to~$K^U$, and hence to $D(2,U)$ because $K$ is of normal type. However, as is easily seen, the identity element is the only element of $D(2,U)$ that leaves all points of $L\le U\cap V$ fixed, a contradiction.\eprf\medskip

Again from statement~(i) of Lemma~\ref{310516a} it follows that the number of groups $H$ such that $\soc(G)\le H\trianglelefteq G$ is at most $\log n$, and so is $|\cH_1|$. For every $H$ and each $k\in H^*\times \id_{G\setminus H}$ one can efficiently test whether $k$ is an automorphism of~$\cX$.  Thus, Lemma~\ref{121116a} immediately implies the following statement.

\crllrl{121116y}
Given a central Cayley scheme $\cX$ over an almost simple group $G$ of order $n$, one can test in time $\poly(n)$ whether the group $\aut(\cX)$ is of  normal type, and (if so) find the principal section of~$\aut(\cX)$ within the same time.
\ecrllr

\section{A majorant for the coset of isomorphisms}\label{061216g}

Throughout this section, we assume that $\cX$ is a central Cayley scheme over an almost simple group~$G$ and $K=\aut(\cX)$. The principal $K$-section of $G$ and the associated partitions are denoted by $U/L$ and $\fU$ and $\fL$, respectively. The equivalence relations corresponding to the partitions $\fU$ and $\fL$, are denoted by $e_\fU$ and $e_\fL$. Let $\varphi$ be an algebraic isomorphism from $\cX$ onto a Cayley scheme $\cX'$ over an almost simple group~$G'$. Assume that
\qtnl{281116b}
\varphi(e_{\fU^{}})=e_{\fU'}\qaq\varphi(e_{\fL^{}})=e_{\fL'},
\eqtn
where in what follows, the group $K'$, the principal section $U'/L'$, the partitions  $\fU'$ and $\fL'$, and the equivalence relations $e_{\fU'}$ and $e_{\fL'}$ are defined for the scheme~$\cX'$ in a similar way.

\lmml{281116a}
In the above notation, $|\fU|=|\fU'|$ and $|\fL|=|\fL'|$. Moreover, the groups $K$ and $K'$ either both of symmetric type, or both of normal type.
\elmm
\proof The first statement follows from assumption~\eqref{281116b}. To prove the second one, we note that by Lemma~\ref{111116e} the group $K$ is of symmetric type if and only if the scheme $\cX$ is isomorphic to the wreath product $\cT_U\wr\cY$, where $\cY$ is the quotient of~$\cX$ modulo the equivalence relation~$e_\fL$. Since algebraic isomorphisms respect wreath products, we are done.\eprf\medskip

For all $Y\in\fU$ and $Y'\in\fU'$, the algebraic isomorphism $\varphi$ induces an algebraic isomorphism
$$
\varphi_{Y^{},Y'}:\cX^{}_{Y^{}}\to\cX'_{Y'},
$$
that takes a relation  $s_{Y}$ to the relation $s'_{Y'}$ for all basis relations $s\subseteq e_\fU$ of the scheme $\cX$, where $s'=\varphi(s)$. It follows that if $\cX'=\cX$ and $\varphi$ is trivial, then $\varphi_{Y,Y}$ is trivial for all $Y\in\fU$.\medskip

For each $Y\in\fU$, set
$$
D_Y=\css
\sym(Y) &\text{if $K$ is of symmetric type},\\
D(2,U^Y)\cap\aut(\cX_Y) &\text{otherwise},\\
\ecss
$$
where, for brevity, $U^Y$ denotes the restriction of the permutation group $U_{right}$ to the set~$Y$. Note that the form of the group $D_Y$ does not depend on $Y\in\fU$, and $D_Y$ contains~$K^Y$ for all~$Y$ (Theorems~\ref{091116a} and~\ref{091116y}). Furthermore, if $Z\in\fU$, then any permutation of~$G^*$ taking $Y$ to $Z$ induces a permutation isomorphism from the group $D_Y$ onto the group $D_Z$.

\lmml{151116a}
$(D_Y)^f=D_{Y^f}$ for all $Y\in\fU$ and all $f\in\iso(\cX,\cX',\varphi)$.
\elmm
\proof Without loss of generality, we may assume that both $K$ and $K'$ are of normal types (Lemma~\ref{281116a}). Let $Y\in\fU$ and  $f\in\iso(\cX,\cX',\varphi)$. Then in view of~\eqref{281116b},
\qtnl{171116a}
\aut(\cX^{}_{Y^{}})^f=\aut(\cX'_{Y'}),
\eqtn
where $Y'=Y^f$. Thus, it suffices to verify that $D(2,V)^f=D(2,V')$, where $V=U^Y$ and $V'=(U')^{Y'}$. However,  by the definition of $\fU$, we have $V=Ug$ for suitable $g\in G$ and
$$
D(2,V)=D(2,U^{g_r}),
$$
where $g_r:U\to V$ is the bijection induced by right multiplication by~$g$. A similar statement holds for $V'$, $U'$, and a suitable bijection $(g')_r$. Therefore, without loss of generality, we may assume that $Y=U$ and $Y'=U'$.\medskip

Recall that  $U$ is an almost simple group with $\soc(U)=L$, and also $U^*\le K^U$. Therefore, $\soc(K^U)=L^*$. The same is true with $K$, $U$, and $L$, replaced by $K'$, $U'$, and $L'$, respectively. Taking into account that $f$ takes $K^U$ to $(K')^{U'}$, we conclude that
$$
(L^*)^f=(L')^*.
$$
This implies that $f$ takes the normalizer of the group $L^*\le\sym(U)$ to the normalizer of the group $(L')^*\le\sym(U')$. However by Theorem~\ref{240916a}, these normalizers are equal to $D(2,U)$ and $D(2,U')$. Thus, $f$ takes the first of these groups to the second, and we are done.\eprf

\dfntnl{171116j}
Denote by $C_\varphi(\cX,\cX')$ the set of all bijections $f:G\to G'$ taking $\fU$ to $\fU'$ and satisfying the following conditions for every $Y\in\fU$:
\qtnl{161116c}
f^Y\in\iso(\cX^{}_{Y^{}},\cX'_{Y^f},\varphi_{Y^{},Y^f})\qaq (D_Y)^f=D_{Y^f}.
\eqtn
\edfntn

Let us find the explicit form of the set $C=C_\varphi(\cX,\cX')$ when $\cX=\cX'$ and $\varphi=\id$. In this case, $C$ is obviously a subgroup of $\sym(G)$ preserving the partition~$\fU$. Condition \eqref{161116c} means that $f^Y$ belongs to the intersection of $\aut(\cX_Y)$ and the normalizer of $D_Y$ in $\sym(Y)$ for all $f\in C_\fU$. This proves the first of the equalities
$$
\prod_{Y\in\fU}D_Y=C_\fU\qaq C^\fU=\sym(\fU);
$$
the second equality follows from the fact that any $g\in G$ induces the permutation isomorphism from $D_{Y^{}}$ onto $D_{Y^g}$ that induces $\varphi_{Y^{},Y^g}$. Thus, the group $C$ is permutation isomorphic to the wreath product $D_U\wr\sym(\fU)$ (in imprimitive action). For arbitrary~$\cX'$ and $\varphi$, an for each $f\in C_\varphi(\cX,\cX')$, we obviously have
\qtnl{161116g}
C_\varphi(\cX,\cX')=C_{\id}(\cX,\cX)f.
\eqtn
Thus if the set $C_\varphi(\cX,\cX')$ is not empty, then it can be given by a generator set of the group $C=C_{\id}(\cX,\cX)$ and the bijection~$f$. In the sense of the following statement, the set $C_\varphi(\cX,\cX')$ can be called a {\it majorant} of $\iso(\cX,\cX',\varphi)$.

\thrml{161116a}
$\iso(\cX,\cX',\varphi)\subseteq C_\varphi(\cX,\cX')$. Moreover, the set
$C_\varphi(\cX,\cX')$ can be found in time $\poly(n)$.
\ethrm
\proof The first statement immediately follows from Lemma~\ref{151116a}. To prove the second one, it suffices to find the set
$$
C_0=\{f_0\in \iso_\psi(\cX^{}_{U^{}},\cX'_{U'}):\ (D_{U^{}})^{f_0}=D_{U'}\},
$$
where $\psi=\varphi_{U^{},U'}$. Indeed, if this set is empty, then obviously so is the majorant $C_\varphi(\cX,\cX')$. On the other hand, if $f_0\in C_0$, then to construct the majorant given by formula~\eqref{161116g} it suffices to find $D_U$ and the bijection $f$ defined as follows:
$$
f^Y=(g_Y)^{-1}f_0g'_{Y'},\qquad Y\in\fU,
$$
where $Y\mapsto Y'$ is an arbitrary bijection from $\fU$ onto $\fU'$ taking $U$ to $U'$, and the bijections $g_Y:U\to Y$ and $g'_{Y'}:U'\to Y'$ are induced by the right multiplications by the elements $g\in G$ and $g'\in G'$ such that $Y=Ug$ and $Y'=U'g'$.\medskip

To find the sets $D_U$ and $C_0$, assume first that $K$ is of symmetric type (recall that this can efficiently checked by Theorem~\ref{121116b}). Then the coherent configurations $\cX^{}_{U^{}}$ and $\cX'_{U'}$ are trivial. Thus, $D_U=\sym(U)$ and for any bijection $f_0:U\to U'$,
$$
C_0=\sym(U)f_0.
$$

Let now $K$ be of normal type. Then $D_U\le D(2,U)$ and $D_{U'}\le D(2,U')$. In particular, $D_U=D(2,U)\cap\aut(\cX_U)$ can be found in time $\poly(n)$. Furthermore, every element $f_0\in C_0$ takes $D_U$ to $D_{U'}$ (Lemma~\ref{151116a}), and induces a permutation group isomorphism from $U_{right}$ onto a group
$$
V'\in\reg(D_{U'},U').
$$
By statement~(ii) of Lemma~\ref{310516a}, the set $\reg(D_{U'},U')$ is of cardinality at most $|U'|^c\le n^c$ for some constant $c>0$, and all its elements can be found by exhaustive search of all $3$-generated subgroups of the group $D_{U'}$. Since for a fixed $V'$, there are at most $|\aut(V')|\le n^c$ distinct elements $f_0\in C_0$ taking $U_{right}$ to  $V'$, one can test in time $\poly(n)$, whether the set $C_0$ is not empty and (if so) find it in the form
$$
C_0=D_Uf_0
$$
with arbitrary $f_0\in C_0$.\eprf

\section{Proof of Theorem~\ref{140516a}}\label{140516w}

\sbsnt{Reduction to Cayley schemes.}
Let $\Gamma$ be a central Cayley graph over an almost simple group $G$,  $\{e_i:\ i\in I\}$ the set of color classes of $\Gamma$, and $K=\aut(\Gamma)$. The principal $K$-section of $G$ and the associated partitions are denoted by $U/L$ and $\fU$ and $\fL$, respectively. Set
$$
\cX=\WL(\Gamma,\{e_\fU,e_\fL\}),
$$
where  $e_\fU$ and $e_\fL$ are the equivalence relations corresponding to the partitions $\fU$ and $\fL$. For other central Cayley graphs~$\Gamma'$, we use similar notation, e.g., $G'$ and $K'$ denote the underlying group and the automorphism group of $\Gamma'$, respectively.

\lmml{thbw}
Given central Cayley graphs $\Gamma$ and $\Gamma'$ over almost simple groups~$G$ and~$G'$, respectively, one can construct in time $\poly(n)$ the Cayley schemes~$\cX$ and $\cX'$ over the same underlying groups and check whether there exists a (unique) algebraic isomorphism $\varphi:\cX\to\cX'$ such that
\qtnl{230516b}
\varphi(e^{}_i)=e'_i\ \,\text{for all}\ i\in I\qquad\text{and}\qquad \varphi(e_{\fU^{}})=e_{\fU'},\  \varphi(e_{\fL^{}})=e_{\fL'},
\eqtn
and (if so) find $\varphi$ within the same time. Moreover, $K=\aut(\cX)$, $K'=\aut(\cX')$, and also
\qtnl{151116u}
\iso(\Gamma,\Gamma')=\iso(\cX,\cX',\varphi).
\eqtn
\elmm
\proof By Theorem~\ref{121116b}, the principal sections and hence the equivalence relations $e_{\fU^{}}$, $e_{\fL^{}}$ and $e_{\fU'}$, $e_{\fL'}$ can be found in time $\poly(n)$. Therefore the first part of the statement immediately follows from Theorem~\ref{131116t}. To prove the second one, we observe that  every $f\in\iso(\Gamma,\Gamma')$ takes the group $K$ to the group $K'$. By the definition of the minimal block (Subsection~\ref{290516a}) this implies that $f$ takes $e_\fL$ to $e_{\fL'}$ and hence takes $e_\fU$ to  $e_{\fU'}$. Thus, the isomorphism~$f$ induces $\varphi$. This means that $f\in\iso(\cX,\cX',\varphi)$ and hence the left-hand side of~\eqref{151116u} is contained in the right-hand side. Since the reverse inclusion is obvious, equality~\eqref{151116u} is completely proved. Next, if $\Gamma'=\Gamma$, then $\varphi=\id$ and  equality~\eqref{151116u} shows that  $K=\aut(\cX)$. Similarly, $K'=\aut(\cX')$.\eprf

\sbsnt{Determining the coset of isomorphisms.}
Denote by $\pi$ the canonical epimorphism from $G$ onto $\ov G=G/L=\fL$. Then $\pi$ induces a map taking the set $\cS$ of basis relations of the Cayley scheme~$\cX$ over $G$ to the set $\ov\cS$ of basis relations of the quotient Cayley scheme  $\ov\cX$ over $\ov G$. In particular, $\pi$ takes $\cS^\cup$ to $\ov\cS^\cup$. Set $\ov\Gamma$ to be the Cayley graph over $\ov G$ with color classes $\pi(e_i)$, $i\in I$. For any set $C$ of the bijections $f:G\to G'$ taking $e_\fU$ to $e_{\fU'}$ and $e_\fL$ to $e_{\fL'}$, we denote by $\ov C$ the set of bijections $\ov f:\ov G\to \ov G'$ induced by~$f\in C$.

\thrml{220516a}
Let $\Gamma$ and $\Gamma'$ be central Cayley graphs over almost simple groups~$G$ and~$G'$, respectively. Assume that the algebraic isomorphism $\varphi$ from Theorem~\ref{thbw} does exist. Then
\qtnl{230516a}
\iso(\Gamma,\Gamma')=\pi^{-1}(\ov C\cap B),
\eqtn
where $B=\iso(\ov\Gamma,\ov\Gamma')$, $C=C_\varphi(\cX,\cX')$, and the right-hand side consists of all $f\in C$ for which $\ov f\in B$.
\ethrm
\proof To prove that that the left-hand side of~\eqref{230516a} is contained
in the right-hand side, let $f\in\iso(\Gamma,\Gamma')$. Then the uniqueness
of the principal sections implies that
$$
(e_{\fU^{}})^f=e_{\fU'}\qaq (e_{\fL^{}})^f=e_{\fL'}.
$$
Therefore, the isomorphism $f$ induces the algebraic isomorphism~$\varphi$. By Theorem~\ref{161116a}, this implies that $f\in C$. Consequently, the induced bijection $\ov f:\ov G\to \ov G'$ belongs to the set~$\ov C$.
Since obviously $\ov f\in B$, we conclude that $f$ belongs to the right-hand side of~\eqref{230516a}, as required.\medskip

Conversely, let $f$ belong to the right-hand side of \eqref{230516a}. By formula~\eqref{151116u} in Lemma~\ref{thbw}, it suffices to verify that $f$ induces $\varphi$. To this end, let $s\in \cS$. Assume first that $s\subseteq e_\fU$. Then $s$ equals the union of $s_Y$, $Y\in\fU$. Therefore since $f\in C$, conditions~\eqref{161116c} are satisfied for all $Y$ and hence
$$
s^f=(\bigcup_{Y\in\fU}s_Y)^f=\bigcup_{Y\in\fU}(s_Y)^f
=\bigcup_{Y\in\fU}\varphi_{Y,Y^f}(s_Y)=\varphi(s).
$$
Now assume that $s$ is outside the equivalence relation $e_\fU$. Let us prove that
\qtnl{300317a}
s=\bigcup_{\substack{X,Y\in\fL,\\(X,Y)\in \ov s}}X\times Y\quad\qaq\quad
s^f=\bigcup_{\substack{X',Y'\in\fL',\\(X',Y')\in \ov s'}}X'\times Y',
\eqtn
where $\ov s=\pi(s)$ and $\ov s'=\pi'(s^f)$. Since $(e_{\fU^{}})^f=e_{\fU'}$, the relation $s^f$ is outside the equivalence relation~$e_{\fU'}$. Therefore, it suffices
to verify the first equality of~\eqref{300317a}; denote the right-hand side of this equality by $t$. Clearly, $s\subseteq t$. Conversely, let $(X,Y)\in t$ for some $(X,Y)\in\ov s$. Since $s$ is outside $e_\fU$, we conclude that $X\times Y$ is a basis relation of the coherent configuration $\inv(K^{}_{\fL^{}})$ (Theorem~\ref{200416b}). Since
$\inv(K^{}_{\fL^{}})\ge \cX$, it follows that this basis relation is contained in a basis relation of $\cX$ which equals $s$, because $(X,Y)\in\ov s$. Therefore, $t\subseteq s$. This completes the proof of~\eqref{300317a} implying $s=\pi^{-1}(\ov s)$ and $s^f=(\pi')^{-1}(\ov s')$.\medskip

On the other hand, the graph isomorphism $\ov f\in B\cap \ov C$ induces an algebraic isomorphism~$\ov\varphi:\ov\cX\to\ov\cX'$ that coincides with the restriction of the algebraic isomorphism~$\varphi$ modulo~$e_{\fL}$. Thus,
$$
\varphi(s)=\varphi(\pi^{-1}(\ov s))=
(\pi')^{-1}(\ov\varphi(\ov s))=(\pi')^{-1}((\ov s)^{\ov f})=
(\pi')^{-1}(\ov s')=s^f,
$$
as required.\eprf

\sbsnt{The algorithm.}
In the algorithm below, the input is given by two central Cayley graphs $\Gamma$ and $\Gamma'$ over almost simple groups $G$ and $G'$, respectively. It is assumed that these groups are
presented by the multiplication tables. The output consists of the set $\iso(\Gamma,\Gamma')$, which is either empty or equals the set $\aut(\Gamma)f$ for some  $f\in\iso(\Gamma,\Gamma')$. Here, the group $\aut(\Gamma)$ is presented by a generating set.
\medskip

\centerline{\bf Central Cayley graph isomorphism test}
\medskip
\noindent{\bf Step 1.}
Find the principal sections of the automorphism groups of the (central Cayley) schemes $\WL(\Gamma)$ and $\WL(\Gamma')$ (Theorem~\ref{121116b}); denote by $\fU$, $\fL$ and $\fU'$, $\fL'$ the associated partitions of $G$ and $G'$, respectively.\vspace{1mm}

\noindent{\bf Step 2.}
Find the schemes $\cX=\WL(\Gamma,\{e_\fU,e_\fL\})$ and $\cX'=\WL(\Gamma',\{e_{\fU'},e_{\fL'}\})$ and the algebraic isomorphism $\varphi$ satisfying condition~\eqref{230516b}; if $\varphi$ does not exist, output $\iso(\Gamma,\Gamma')=\varnothing$.
\vspace{1mm}

\noindent{\bf Step 3.}
Find the set $C=C_\varphi(\cX,\cX')$ (Theorem~\ref{161116a}).
\vspace{1mm}

\noindent{\bf Step 4.} Using the graph isomorphism and coset intersection  algorithms from \cite{BKL}, find the set $B=\iso(\ov\Gamma,\ov\Gamma')$ and then the set $B'=B\cap\ov C$.\vspace{1mm}

\noindent{\bf Step 5.} Output $\iso(\Gamma,\Gamma')=\pi^{-1}(B')$.\eprf\medskip

To complete the proof of Theorem~\ref{140516a}, we show that the above algorithm correctly finds the set $\iso(\Gamma,\Gamma')$ in time $\poly(n)$. Note that every graph  isomorphism $f\in\iso(\Gamma,\Gamma')$ induces an algebraic isomorphism $\varphi$ satisfying condition~\eqref{230516b}.
Therefore, the output at Step~2 is correct. Thus, the correctness of the output at
Step~5 and hence of the algorithm immediately follows from Theorem~\ref{220516a}.\medskip

To estimate the running time, we note that all the steps except for Step~4 run in polynomial time (Theorem~\ref{121116b}, Lemma~\ref{thbw}, and Theorem~\ref{161116a}). Furthermore, the graph isomorphism and coset intersection  algorithms from \cite{BKL} are applied at Step~4 to graphs with $m=|\fL|$ vertices and to the cosets contained in $\sym(m)$, respectively. 
Each of these algorithms runs in time at most $\exp((\log m)^c\sqrt{m})$.
Since $m\le \log n$,  the complexity of  this step does not exceed
$$
\exp((\log m)^c\sqrt{m})\le \exp((\log\log n)^c(\log{n})^{1/2})
\le  \exp(\log{n})\le n
$$
for sufficiently large $n$ and a suitable constant $c>0$. Thus, the running time of the algorithm is polynomial in~$n$, as required.\eprf

\end{document}